\documentclass[11pt]{article}
\usepackage{amssymb} 

\def\a{\alpha}

\def\f{\varphi}
\def\g{\gamma}

\def\i{\iota}

\def\l{\lambda}
\def\o{\omega}

\def\r{\rho}
\def\s{\sigma}

\def\u{\upsilon}

\chardef\tempcat=\the\catcode`\@
\catcode`\@=11
\def\cyracc{\def\u##1{\if \i##1\accent"24 i
    \else \accent"24 ##1\fi }}
\newfam\cyrfam
\font\tencyr=wncyr10
\def\cyr{\fam\cyrfam\tencyr\cyracc}

\font\cmsslll=cmss10 at 14 pt

\DeclareFontFamily{OT1}{msb}{}{}
\DeclareFontShape{OT1}{msb}{m}{n}
 {  <5> <6> <7> <8> <9> <10> gen * msbm
      <10.95><12><14.4><17.28><20.74><24.88>msbm10}{}
\DeclareMathAlphabet{\bubble}{OT1}{msb}{m}{n}

\def\bR{{\bubble R}}

\def\bC{{\bubble C}}
\def\bH{{\bubble H}}

\newfont{\goth}{eufm10 scaled \magstep1}
\def\ga{\mbox{\goth a}}
\def\gb{\mbox{\goth b}}
\def\gc{\mbox{\goth c}}

\def\gg{\mbox{\goth g}}
\def\gh{{\mbox{\goth h}}}
\def\gk{\mbox{\goth k}}
\def\gl{\mbox{\goth l}}
\def\gm{\mbox{\goth m}}

\def\gp{\mbox{\goth p}}

\def\gr{\mbox{\goth r}}
\def\gs{\mbox{\goth s}}

\newfont{\mcal}{eusm10 scaled \magstep1}

\def\square{\kern1pt\vbox
            {\hrule height 0.6pt\hbox{\vrule width 0.6pt\hskip 3pt
 \vbox{\vskip 6pt}\hskip 3pt\vrule width 0.6pt}\hrule height 0.6pt}\kern1pt}

\def\ra{\rightarrow}

\def\Hol{\mathrm{Hol\;}}

\newtheorem{Th}{Theorem}
\newtheorem{Prop}{Proposition} 
\newtheorem{Cor}{Corollary} 
\newtheorem{Lem}{Lemma}
\newtheorem{Def}{Definition} 

\def\bt{\begin{Th}}
\def\et{\end{Th}}
\def\bp{\begin{Prop}}
\def\ep{\end{Prop}}
\def\bc{\begin{Cor}}
\def\ec{\end{Cor}}
\def\bl{\begin{Lem}}
\def\el{\end{Lem}}
\def\bd{\begin{Def}}
\def\ed{\end{Def}}

\def\pf{\noindent{\it Proof:\ }}
\def\qed{\hfill\square}
\def\n{\nabla}
 \def\ot{\otimes}

\def\be{\begin{equation}}
\def\ee{\end{equation}}
\def\la#1{\label{#1}} 
 
\def\arr{\begin{array}{rlll}}
\def\ea{\end{array}}
\def\bea{\begin{eqnarray}}
\def\eea{\end{eqnarray}}  
\def\bean{\begin{eqnarray*}}
\def\eean{\end{eqnarray*}}  

\catcode`@=11
\@addtoreset{equation}{section}
\catcode`@=12   

\def\3s{3-Sasakian manifold} 
\def\hk{hyper-K\"ahler\, }
\def\hkss{hyper-K\"ahler symmetric space\, }

\begin{document}
\begin{center}
{\LARGE Classification of indefinite hyper-K\"ahler symmetric spaces}
\vskip 1.0 true cm
{\cmsslll
Dmitri V. Alekseevsky$\,^{1}$, \ Vicente Cort\'es$\,^{2}$}
\vskip 0.8 true cm
{\small daleksee@mpim-bonn.mpg.de}\ ,\ 
{\small vicente@math.uni-bonn.de}\ ,\\

{\cyr $^{1}$  Centr ``Sofus Li'', Gen.\ Antonova 2 - 99, 117279 Moskva} 
{\small and 
MPI f\"ur Mathematik, Bonn}\\[2pt] 
{\small $^2$ Mathematisches Institut der Universit\"at Bonn,
Beringstr. 1, D-53115 Bonn}
\end{center}
\begin{quote} 
\centerline{\bf Abstract}
We classify indefinite  simply connected hyper-K\"ahler
symmetric spaces. Any such space without flat factor 
has commutative holonomy group and signature $(4m,4m)$. 
We establish a  natural 1-1 correspondence between simply connected 
hyper-K\"ahler symmetric spaces of dimension $8m$ and orbits of the group 
GL$(m,\bH )$ on the space $(S^4\bC^n)^{\tau}$ of homogeneous quartic 
polynomials $S$ in $n = 2m$ complex variables satisfying the reality
condition $S = \tau S$, where $\tau$ is the real structure induced
by the quaternionic structure of $\bC^{2m} = \bH^m$. We define and
classify also complex hyper-K\"ahler
symmetric spaces. Such spaces without flat factor exist in any 
(complex) dimension divisible by $4$. 
\end{quote}
\def\large{\bf}
\section{Introduction}
We recall that a pseudo-Riemannian manifold $(M,g)$ is called a
symmetric space if any point $x \in M$ is an isolated fixed point
of an involutive isometry $s_x$ (called central symmetry with centre
$x$). Since the product of two central symmetries $s_x$ and $s_y$
with sufficiently close centres  is a shift along the geodesic $(xy)$,
the group generated by central symmetries acts transitively on $M$
and one can identify $M$ with the quotient $M = G/K$ where $G$ is the
connected component of the isometry group ${\rm Isom}(M,g)$ and $K$
is the stabilizer of a point $o \in M$.

 A symmetric space $(M =G/K, g)$ is called K\"ahler (respectively, 
hyper-K\"ahler) if its holonomy group Hol$(M,g)$ is a subgroup of
the pseudo-unitary group ${\rm U}(p,q)$ (respectively, of the pseudo-symplectic
group ${\rm Sp}(p,q) \subset {\rm SU}(2p,2q)$). Any hyper-K\"ahler 
symmetric space is in particular a homogeneous hypercomplex manifold. 
Homogeneous hypercomplex manifolds 
of compact Lie groups were constructed by Ph.\ Spindel, A.\ Sevrin,
W.\ Troost, A.\ Van Proeyen \cite{SSTVP} and by D.\ Joyce \cite{J} and  
homogeneous hypercomplex structures on solvable  Lie groups by 
M.L.\ Barberis and I.\ Dotti-Miatello \cite{BD}. 

 The classification of simply connected symmetric spaces 
reduces to  the classification of involutive automorphisms $\sigma$ 
of a Lie algebra $\gg$, such that the adjoint representation
${\rm ad}_{\gk}|\gm$ preserves a pseudo-Euclidean
scalar product $g$, where 
$$ \gg = \gk + \gm \, , \quad \sigma |\gk = 1, \quad \sigma |\gm = -1 ,$$ 
is the  eigenspace decomposition of the involution $\sigma$.
Note that the eigenspace decomposition of an involutive automorphism
is characterized by the conditions
$$ [\gk, \gk ] \subset \gk ,\quad 
[\gk, \gm] \subset \gm , \quad [\gm, \gm]\subset
\gk \, .$$
Such a decomposition is called a symmetric decomposition.

In fact, for any pseudo-Riemannian symmetric space $M = G/K$
the conjugation with respect to the central symmetry $s_o$ 
with centre $o = eK$ is an involutive automorphism  of the Lie group
$G$, which induces an involutive automorphism $\sigma$ 
of its Lie algebra $\gg$. The pseudo-Riemannian metric of $M$ induces
a $\gk$-invariant scalar product on $\gm \cong T_oM$, where $\gg = \gk + \gm$
is the symmetric decomposition defined by $\sigma$. 
Conversely, a symmetric decomposition
$\gg = \gk + \gm$ together with a $\gk$-invariant scalar product on $\gm$ 
determines a pseudo-Riemannian symmetric space 
$M = G/K$ where $G$ is the 
simply connected Lie group with the Lie algebra $\gg$, $K$ is the connected
(and closed) subgroup of $G$ generated by $\gk$, the pseudo-Riemannian
metric on $M$ is defined by $g$ and the central symmetry is defined by
the involutive automorphism $\sigma$ associated to the symmetric 
decomposition. 

Naturally identifying the space $\gm$ with the tangent space
$T_oM$, the isotropy group is identified with
${\rm Ad}_{K}|\gm$ and the holonomy algebra is identified with
${\rm ad}_{\gh}$ where $ \gh = [\gm, \gm]$.
If one assumes that  the holonomy algebra is irreducible then one can prove
that the Lie algebra $\gg$ is semisimple. Hence the 
classification of pseudo-Riemannian symmetric spaces with irreducible
holonomy reduces to the classification of involutive authomorphisms of
semisimple Lie algebras. Such a classification was obtained by M.\ Berger
\cite{B1,B2} and A.\ Fedenko \cite{F}. It includes the classification of 
Riemannian symmetric spaces (obtained earlier by E.\ Cartan), 
since according to de
Rham's theorem any simply
connected complete Riemannian manifold is a direct product of
Riemannian manifolds with irreducible holonomy algebra and a Euclidean
space.   

A classification of pseudo-Riemannian symmetric spaces with non
completely reducible holonomy is known only for signature
$(1,n)$ (Cahen-Wallach \cite{C-W}) and for signature $(2,n)$ 
under the assumption that the holonomy
group is solvable (Cahen-Parker \cite{C-P}). The classification problem for
arbitrary signature 
looks very complicated and includes, for example, the classification
of Lie algebras which admit a nondegenerate ad-invariant symmetric bilinear
form. Some inductive construction of solvable Lie algebras with such
a form was given by V.\ Kac. 
  
In this paper we give a classification of pseudo-Riemannian hyper-K\"ahler
symmetric spaces.  In particular, we prove that any simply connected  
hyper-K\"ahler symmetric space $M$ 
has signature (4m,4m) and  its holonomy group is commutative. 
The main result is the following, see Theorem \ref{realmainThm}.

 Let $(E,\omega ,j)$  be a complex symplectic vector space of
dimension $4m$ with a quaternionic structure $j$ such that 
$\o (j x ,j y ) = \overline{\o (x,y)}$ for all $x,y\in E$ 
and $E = E_+ \oplus E_-$ a $j$-invariant 
Lagrangian decomposition. Such a decomposition exists if and only if 
the Hermitian form $\gamma = \o (\cdot , j\cdot )$  has real signature
$(4m,4m)$.  We denote by $\tau$ the real structure in $S^{2r}E$ defined by
$\tau ( e_1 e_2 \dots e_{2r}) := j(e_1)j(e_2)\dots j(e_{2r})$, $e_i \in E$. 
Then any element $S \in (S^4E_+)^{\tau}$
defines a hyper-K\"ahler symmetric space $M_S$ which is associated
with the symmetric decomposition
$$ \gg = \gh + \gm \, ,$$ 
where $\gm = (\bC^2 \otimes E)^{\rho}$
is the fixed point set of the real structure 
$ \rho $ on $\gm$ given by 
$\rho(h\ot e) = j_Hh \ot je$, where $j_H$ is the standard quaternionic
structure on $\bC^2 = \bH$, 
$\gh ={\rm span}\{S_{e e'}|\, e, e' \in E \}^{\tau}\subset
sp(E)^{\tau} \cong {\rm sp}(m,m) $ with the natural action on 
$\gm \subset \bC^2 \ot E $
and the Lie bracket $\gm \wedge \gm \ra \gh $
is given by
$$ [h \ot e, \, h'\ot e'] = \omega_H (h, h')S_{e, e'} .$$ 

 Moreover, we establish a natural 1-1 correspondence between
simply connected hyper-K\"ahler symmetric spaces 
 (up to isomorphism) and
orbits of the group $GL(m,\bH)$ in   
$(S^4E_+)^{\tau}$.

 We define also the notion of complex hyper-K\"ahler symmetric space
as a complex manifold $(M,g)$ of complex dimension $4n$ with holomorphic
metric $g$ such that for any point $x \in M$ there is a holomorphic  
central symmetry $s_x$ with centre $x$ and which 
has holonomy group ${\rm Hol}(M,g) \subset {\rm Sp}(n, \bC)$ 
(${\rm Sp}(n, \bC)\hookrightarrow {\rm Sp}(n,\bC ) \times {\rm Sp}(n,\bC ) 
\subset {\rm O}(4n, \bC)$ is diagonally embedded) and give a
classification of such spaces.   We establish a natural 1-1 correspondence
between simply connected complex hyper-K\"ahler symmetric spaces 
and homogeneous polynomials of degree 4 in the vector space
$\bC^{2n}$ considered up to linear transformations from $GL(2n, \bC)$.

\section{Symmetric spaces} 
\subsection{Basic facts about pseudo-Riemannian symmetric spaces}
A {\bf pseudo-Riemannian symmetric space} is a pseudo-Riemannian manifold 
$(M,g)$ such that any point is an isolated fixed point of an isometric 
involution. Such a pseudo-Riemannian manifold admits a transitive Lie group of
isometries $L$ and can be identified with $L/L_o$, where $L_o$ is the 
stabilizer of a point $o$.   More precisely, any simply connected
pseudo-Riemannian symmetric 
space $M = G/K$ is associated with a symmetric decomposition 
\be \gg = \gk + \gm \, ,\quad [\gk ,\gk ] \subset \gk\, ,\quad 
[\gk ,\gm ] \subset \gm \, , \quad 
[\gm ,\gm ] \subset \gk \label{symdecomp}\ee  
of the Lie algebra $\gg = Lie\, G$ together with an ${\rm Ad}_K$-invariant
pseudo-Euclidean scalar product on $\gm$. We will assume that  G acts
almost effectively on $M$, i.e.\ $\gk$ does not contain any nontrivial ideal
of $\gg$, that
$M$ and $G$ are simply connected and that $K$ is connected.  
Then, under the natural identification of the tangent space 
$T_oM$ 
at the canonical base point $o = eK$ with $\gm$, the holonomy group 
$\Hol \subset $ Ad$_K|\gm$.  We will denote by $\gh$ the 
holonomy Lie algebra. Since the 
isotropy representation
is faithfull it is identified with the subalgebra 
$\gh = [\gm , \gm ] := {\rm span}\{ [x,y]|\, x,y \in \gm\} \subset \gk$.  Recall that the curvature
tensor $R$ of a symmetric space $M$ at $o$ is $\gh$-invariant and determines 
the Lie bracket in the ideal $\gh + \gm \subset \gg$ as follows:
\[ \gh = R(\gm , \gm ) := {\rm span} \{ R(x,y)| x,y \in \gm\} \quad \mbox{and} 
\quad [x,y] = -R(x,y)\, ,\quad x,y \in \gm \, .\] 
The following result is well known: 
\bp \label{fullisotropyProp} 
The full Lie 
algebra of Killing fields of a symmetric space has the form
\[ {\rm isom}(M) = \tilde{\gh} + \gm \, ,\] 
where the full isotropy  subalgebra is given by  
\be \tilde{\gh} = {\rm aut}(R) = \{ A \in {\rm so}(\gm )| 
A \cdot R = [A, R(\cdot ,\cdot )] - R(A \cdot ,\cdot ) 
 - R(\cdot ,A\cdot ) = 0\}\, . \label{fullisotropyEq} \ee 
\ep

\subsection{Symmetric spaces of semisimple Lie groups} 
We will prove that in the case when $(M = G/K,g)$ is a pseudo-Riemannian
symmetric space of a (connected) semisimple Lie group $G$ then 
$G$ is the maximal connected Lie group of isometries of $M$.  

\bp \label{ssfullisomProp} Let $(M = G/K,g)$ be a pseudo-Riemannian
symmetric space associated with a symmetric decomposition $\gg = \gk +\gm$.
If $G$ is semisimple and almost effective then 
\begin{enumerate}
\item[(i)] the restriction of the Cartan-Killing form $B$ of $\gg$ to 
$\gk$ is nondegenerate and hence $\gk$ is  a reductive subalgebra of $\gg$
and $\gg = \gk +\gm$ is a $B$-orthogonal decomposition, 
\item[(ii)] $\gk = [\gm ,\gm ]$  and 
\item[(iii)] $\gg = {\rm isom} (M,g)$ is the Lie algebra of the full isometry
group of $M$. 
\end{enumerate}
\ep 

\pf
For (i) see \cite{O-V} Ch.\ 3 Proposition 3.6.\\   
(ii) It is clear that $\bar{\gg} = [\gm ,\gm ] + \gm$ is an ideal of $\gg$. 
The $B$-orthogonal complement $\ga := \bar{\gg}^{\perp}
\subset \gk$ is a complementary ideal of $\gg$. Since $[ \ga , \gm ] = 0$ 
the Lie algebra $\ga$ acts trivially on $M$. From the effectivity 
of $\gg$ we conclude that $\ga = 0$.\\
(iii)  By Proposition \ref{fullisotropyProp}, 
$\tilde{\gg} = {\rm isom} (M,g) = 
\tilde{\gh} + \gm$, where $\tilde{\gh} = {\rm aut}(R) = 
\{ A \in {\rm so}(\gm )| A \cdot R = 0\}$. Now 
$\tilde{\gh}$ preserves $\gm$ and by the identity
$A \cdot R = [A, R(\cdot ,\cdot )] - R(A \cdot ,\cdot ) 
 - R(\cdot ,A\cdot )$ it also normalizes $\gk$. This shows that 
$\tilde{\gh}$ normalizes $\gg$ and hence $\gg \subset \tilde{\gg}$ is an 
ideal. Since $\gg$ is semisimple there exists a $\gg$-invariant complement
$\gb$ in $\tilde{\gg}$. Note that $[\gg ,\gb ] \subset \gg \cap \gb = 0$.
We can decompose any  $X\in \gb$ as $X = Y + Z$, 
where $Y \in \tilde{\gh}$ and $Z \in \gm$. From $[\gg ,\gb ] = 0$ it follows
that $[\gg , Y] = [\gg, Z] = 0$ and in particular $[\gm , Y] = 0$. This 
implies that $Y = 0$ and $X = Z \in \gb \cap \gm = 0$. This shows that
$\gb = 0$ proving (iii). 
\qed

We recall that  a pseudo-Riemannian Hermitian 
symmetric space is pseudo-Riemannian symmetric space
$(M = G/K,g)$ together with an invariant (and hence parallel) $g$-orthogonal 
complex structure $J$. 
\bp Let $(M = G/K,g,J)$ be a pseudo-Riemannian Hermitian 
symmetric space of a semisimple and almost effective Lie group $G$. 
Then the Ricci curvature of $M$ is not zero. \label{Ricneq0} 
\ep   

\pf From Proposition \ref{ssfullisomProp} it follows that $\gg = 
{\rm isom} (M,g) = \tilde{\gh} + m$, where  $\tilde{\gh} = \gk = [\gm ,\gm ]$.
It is well known that the curvature tensor $R$ of any 
pseudo-K\"ahler manifold (and in particular of any pseudo-Riemannian Hermitian
symmetric space) is invariant under the operator $J$. This shows that 
$J \in \tilde{\gh} = {\rm aut}(R) =  [\gm ,\gm ] = \gh$ (holonomy Lie 
algebra), which implies that the holonomy Lie algebra is not a subalgebra of 
${\rm su}(\gm ) \cong {\rm su}(p,q)$. Hence $M$ is not Ricci-flat. 
In fact, we can write
$J = \sum {\rm ad} [X_i,Y_i]$, for $X_i, Y_i \in \gm$. Then using the 
formulas ${\rm Ric}(X,JY) = {\rm tr} JR(X,Y)$ for the Ricci curvature
of a pseudo-K\"ahler manifold and $R(X,Y) = - ad_{[X,Y]}|\gm$ for the 
curvature of a symmetric space we calculate:
\[ \sum {\rm Ric}(X_i,JY_i) = \sum {\rm tr} JR(X_i,Y_i) = - \sum 
{\rm tr} J{\rm ad} [X_i,Y_i] = -{\rm tr} J^2 \neq 0\, .\]
\qed         

\section{Structure of hyper-K\"ahler symmetric spaces} 
\subsection{Definitions} 
A (possibly indefinite) {\bf hyper-K\"ahler manifold} is a 
pseudo-Riemannian manifold $(M^{4n},g)$ of signature $(4k,4l)$ together with 
a {\bf compatible} hypercomplex structure, i.e.\ three $g$-orthogonal parallel 
complex structures $(J_1,J_2,J_3 = J_1J_2)$. This means that the holonomy 
group $\Hol \subset {\rm Sp}(k,l)$.  Two \hk manifolds $(M,g,J_{\alpha})$  
($\alpha  = 1,2,3$)
and $(M',g',J_{\alpha}')$ are called {\bf isomorphic} if there exists a
triholomorphic isometry $\varphi :M \rightarrow M'$, i.e.\ 
$\varphi^*J_{\alpha}' = J_{\alpha}$ and $\varphi^*g'= g$.

A {\bf hyper-K\"ahler symmetric space} is a pseudo-Riemannian symmetric space
$(M = G/K,g)$ together with an invariant  compatible   
hypercomplex structure.  Consider now a simply connected 
hyper-K\"ahler symmetric space
$(M = G/K,g,J_{\alpha})$. Without restriction of generality we will assume 
that $G$ acts almost effectively. 
$M$ being hyper-K\"ahler is equivalent to  
Ad$_K|\gm \subset {\rm Sp}(k,l)$, or, since $K$ is connected, 
to ad$_{\gk}|\gm \subset {\rm sp}(k,l)$.  This condition 
means that $\gk$ commutes with the Lie algebra
$Q = {\rm sp}(1) \subset {\rm so}(\gm ) = {\rm so}(4k,4l)$ spanned by three 
anticommuting 
complex structures $J_1,J_2,J_3$. 

\subsection{Existence of a transitive solvable group of isometries and
solvability of the holonomy} 
In this subsection we prove that any simply connected
hyper-K\"ahler symmetric space
$(M,g,J_\a )$ admits a transitive solvable Lie group $G \subset 
{\rm Aut}(g,J_\a )$ of automorphisms and has solvable holonomy 
group. 

\bp \label{solvableProp} 
Let $(M = G/K,g,J_\a)$ be a simply connected hyper-K\"ahler symmetric space
and $A = {\rm Aut}_0(g,Q) \supset {\rm Aut}_0(g,J_\a ) \supset G$ 
the connected group of isometries which
preserve the quaternionic structure $Q = {\rm span}\{ J_\a \}$. 
Then 
\begin{enumerate}
\item[(i)] the stabilizer $A_o$ of a point $o\in M$ contains a 
maximal semisimple subgroup of $A$,
\item[(ii)] the radical $R$ of $A$ acts transitively and triholomorphically 
on $M$ and 
\item[(iii)] the holonomy group of $M$ is solvable. 
\end{enumerate}
\ep

\pf
We consider the  quaternionic K\"ahler symmetric space $(M = A/A_o ,g,Q)$. The 
Lie algebra $\ga_o$ of the stabilizer is given by
\[ \ga_o = {\rm aut}(R,Q) = \{ A\in {\rm so}(T_oM )| 
A \cdot R = 0, [A,Q] \subset Q\} \, .\]
Since the curvature tensor of a  quaternionic K\"ahler manifold is
invariant under the quaternionic structure $Q$ we conclude that
$Q \subset \ga_o$ and $\ga_o =  Q\oplus Z_{\ga} (Q)$, where $Z_{\ga} (Q)$ 
denotes the centralizer of $Q$ in $\ga$.  Since 
$Q \cong {\rm sp}(1)$ is simple, we may choose a
Levi-Malcev decomposition $\ga = \gs + \gr$ such that the Levi subalgebra
$\gs \supset Q$. We put $\gm_r := [Q,\gr ]$ and denote by
$\gm_s$ a $Q  \oplus Z_{\gs} (Q)$-invariant complement of $Q$ in 
$[Q,\gs ]$. The stabilizer
has the decomposition $\ga_o = Q \oplus (Z_{\gs} (Q) + Z_{\gr} (Q))$. 

\bl 
The complement $\gm = \gm_s + \gm_r$ to $\ga_o$ in $\ga$ 
is  $\ga_o$-invariant and the decomposition
\[ \ga = \ga_o + \gm \]
is a symmetric decomposition.
\el 

\pf
It is clear that $\gm_r$ is $\ga_o$-invariant
and $\gm_s$ is invariant under $Q \oplus Z_{\gs} (Q)$ by construction. 
It remains to check that $[Z_{\gr} (Q),\gm_s] \subset \gm$. 
Since $\gm_s = [Q,\gm_s]$, we have 
\[ [Z_{\gr} (Q),\gm_s] = [Z_{\gr} (Q),[Q,\gm_s]] = [Q,[Z_{\gr} (Q),\gm_s]] \subset [Q,\gr ] = \gm_r \subset \gm \, .\] 
This shows that  $\ga = \ga_o + \gm$ is an $\ga_o$-invariant decomposition.
We denote by $\ga = \ga_o + \gp$ a symmetric decomposition. Any other 
$\ga_o$-invariant decomposition is of the form
$\ga = \ga_o + \gp_{\f}$, where $\f : \gp \rightarrow \ga_o$ is an 
$\ga_o$-equivariant map and $\gp_{\f} = \{ X + \f (X) | X \in \gp \}$. 
If such non-zero equivariant map $\f$ exists then $\gp$ and $\ga_o$ contain
non-trivial isomorphic $Q$-submodules. Since $\gp$ is a sum of 4-dimensional
irreducible $Q$-modules and $\ga_o$ is the sum of the 3-dimensional
irreducible $Q$-module $Q$ and the trivial complementary $Q$-module
$Z_{\ga_o}(Q)$, we infer that there exists a unique $\ga_o$-invariant 
decomposition, which coincides with the symmetric decomposition 
$\ga = \ga_o + \gp$. 
\qed 

To prove (i) we have to check that $\gm_s = 0$. We note that by the 
previous lemma $\gs = (Q \oplus Z_{\gs} (Q)) + \gm_s$ is a symmetric
decomposition of the semisimple Lie algebra $\gs$. Since  $[\gm_s , \gm_s] 
\subset Z_{\gs} (Q)$ it defines a hyper-K\"ahler symmetric space
$N = L/L_o$, where $L$ is the simply connected semisimple Lie group
with Lie algebra $\gl =  Z_{\gs} (Q) + \gm_s$ and $L_o$ is the Lie subgroup 
generated by the subalgebra $Z_{\gs} (Q) \subset \gl$. Since $N$ is 
in particular a Ricci-flat pseudo-Riemannian Hermitian 
symmetric space, from Proposition \ref{Ricneq0} we obtain that $N$ is reduced
to point. Therefore $\gm_s = 0$. This proves (i) and (ii). 
Finally, since the holonomy Lie algebra $\gh$ is identified with 
$\gh = [\gm , \gm ] = [\gm_r , \gm_r ] \subset \gr$ it is solvable
as subalgebra of the solvable Lie algebra $\gr$. 
\qed

\subsection{Hyper-K\"ahler symmetric spaces and second prolongation of 
symplectic Lie algebras}
Let $(M = G/K,g,J_{\alpha})$ be a simply connected 
hyper-K\"ahler symmetric space associated with a symmetric 
decomposition (\ref{symdecomp}). Without restriction of generality 
we will assume that $G$ acts almost effectively and 
that $\gk = [\gm , \gm ] = \gh$ (holonomy Lie algebra).
The complexification $\gm^{\bC}$ as $\gh^{\bC}$-module can be written as 
$\gm^{\bC} = H\ot E$, such that $\gh^{\bC}\subset {\rm Id} \ot 
{\rm sp}(E) \cong {\rm sp}(E)$, where 
$H = {\bC}^2$ and $E = {\bC}^{2n}$ are 
complex symplectic vector spaces with symplectic form $\o_H$ and 
$\o_E$, respectively, such that   
$g^{\bC} = \o_H \ot \o_E$ is the complex bilinear metric on $\gm^{\bC}$
induced by $g$. Note that the symplectic forms are unique up to 
the transformation $\o_H \mapsto \l \o_H$, $\o_E \mapsto \l^{-1} \o_E$, 
$\l \in \bC^*$. 
We have also quaternionic structures 
$j_H$ and $j_E$ on $H$ and $E$,  such that 
$\o_H (j_H x ,j_H y) = \overline{\o_H (x,y)}$ for all $x,y\in H$ and 
$\o_E (j_E x ,j_E y ) = \overline{\o_E(x,y)}$ for all $x,y\in E$, where 
the bar denotes complex conjugation. This implies that $\gamma_H :=
\o_H (\cdot , j_H \cdot )$ and $\gamma_E :=
\o_E (\cdot , j_E \cdot )$ are Hermitian forms on $H$ and $E$. For fixed 
$\o_H$ and $\o_E$ the quaternionic structures 
$j_H$ and $j_E$ are uniquely determined if we require that 
$\gamma_H$ is positive definite and that $\r = j_H \ot j_E$ is 
the real structure on 
$\gm^{\bC}$, i.e.\
the complex conjugation with respect to $\gm$. 
The metric $g^{\bC}$ and the Hermitian form $g^{\bC}(\cdot , \r \cdot ) = 
\gamma_H \ot \gamma_E$ 
restrict to a real valued scalar product $g$ of some signature $(4k,4l)$ on 
$\gm = (H\ot E)^\r$, where $(2k,2l)$ is the (real) signature of the
Hermitian form $\gamma_E = \o_E (\cdot , j_E\cdot )$. Note that 
for the holonomy algebra we have the inclusion   
\[ \gh = {\rm Id} \ot (\gh^{\bC})^{j_E} \hookrightarrow  
{\rm sp}(E)^{j_E} =  \{ A \in  {\rm sp}(E)| [A,j_E] = 0\} = 
{\rm aut}(E,\o_E, j_E ) \cong 
{\rm aut}(\gm , g , J_{\alpha}) \cong {\rm sp}(k,l) \, .\] 
Using the symplectic 
forms we identify
$H = H^*$ and $E = E^*$. Then the symplectic Lie algebras are identified
with symmetric tensors as follows:
\[ {\rm sp}(H) = S^2H\, ,\quad {\rm sp}(E) = S^2E\, .\]

Since the curvature tensor $R$ of any \hk manifold $M^{4n}$ can be identified
at a point $p\in M$ with an element $R \in S^2 {\rm sp}(k,l)$ it is invariant
under the Lie algebra ${\rm sp}(1) ={\rm span}\{ J_1, J_2, J_3\}$.  
Let $M = G/K$ be a \hkss as above.  By Proposition~\ref{fullisotropyProp}
we can extend the Lie algebra $\gg = \gk + \gm = \gh + \gm$ to a Lie algebra 
\[ \tilde{\gg} = {\rm sp}(1) + \gh + \gm \] 
of Killing vector fields such that $[{\rm sp}(1),\gh ] = 0$. 
In the $H \ot E$-formalism the Lie algebra ${\rm sp}(1)$ is identified 
with ${\rm sp}(H)^{j_H} \ot {\rm Id} \subset {\rm so}(\gm )$. 

\bl \label{lemma1} 
Denote by $\tilde{\gg}^{\bC} = {\rm sp}(1,\bC ) + \gh^{\bC } + \gm^{\bC}$
the complexification of the Lie algebra $\tilde{\gg}$. Then the Lie bracket
$[\cdot ,\cdot ] : \wedge^2\gm^{\bC} \rightarrow \gh^{\bC }$ can be written
as 
\be 
[h\ot e, h'\ot e'] = \o_H(h,h')S_{e,e'}\, ,\label{bracketEq}
\ee
where $S \in (\gh^{\bC})^{(2)} := \gh^{\bC}\ot S^2E^* \cap E\ot S^3E^*
=\gh^{\bC}\ot \gh^{\bC} \cap S^4E$. Moreover $S$ is 
${\rm sp}(1,\bC ) \oplus \gh^{\bC }$-invariant and satisfies the following 
reality condition: $[S_{j_Ee,e'}-S_{e,j_Ee'},j_E] = 0$. 
\el 
 
\pf 
The Lie bracket $[\cdot ,\cdot ] : \wedge^2\gm^{\bC} \rightarrow \gh^{\bC }$
is an ${\rm sp}(1,\bC ) \oplus \gh^{\bC }$-equivariant map, due to the Jacobi 
identity. We decompose the ${\rm sp}(H)\oplus {\rm sp}(E)$-module
$\wedge^2\gm^{\bC}$:
\[ \wedge^2\gm^{\bC} =  \wedge^2(H\ot E) = \wedge^2H \ot S^2E \oplus
S^2H \ot \wedge^2E = \o_H\ot S^2E \oplus
S^2H \ot \wedge^2E\, .\]
Since $\gh^{\bC} \subset S^2E$ the Lie bracket defines an 
${\rm sp}(1,\bC ) \oplus \gh^{\bC }$-invariant element of
the space $\o_H\ot S^2E \ot S^2E \oplus
S^2H \ot \wedge^2E\ot S^2E$. The second summand has no 
nontrivial ${\rm sp}(1,\bC )$-invariant elements.  Hence 
the bracket is of the form (\ref{bracketEq}), where 
$S\in S^2E^*\ot \gh^{\bC } \subset S^2E \ot S^2E$. The Jacobi identity
reads:
\[ 0 = [h\ot e, [h'\ot e',h''\ot e'']] - [[h\ot e,h'\ot e'],h''\ot e''] - 
[h'\ot e',[h\ot e,h''\ot e'']] \] 
\[ = 
- \o_H(h',h'')h\ot S_{e',e''}e - \o_H(h,h')h''\ot S_{e,e'}e'' 
+ \o_H(h,h'')h'\ot S_{e,e''}e'\, .\]
Since $\dim H = 2$ we may assume that $h,h'=h''$ is a symplectic basis, i.e.\
$\o_H(h,h')=1$, and the equation implies: $S_{e,e''}e' = S_{e,e'}e''$, i.e.\ 
$S\in (\gh^{\bC})^{(2)}$. The Lie bracket of two real elements 
$h\ot e + j_Hh\ot j_Ee$ and $h\ot e' + j_Hh\ot j_Ee'\in \gm \subset 
\gm^{\bC}$ is an element 
of $\gh$. This gives:
\[ [h\ot e + j_Hh\ot j_Ee,h\ot e' + j_Hh\ot j_Ee'] =
\o_H(h,j_Hh)(S_{e,j_Ee'}-  S_{j_Ee,e'})  \in  \gh  \, .\]
{}From the fact that the Hermitian form $\g_H = \o_H (\cdot , j_H\cdot )$ is 
positive definite  it 
follows that $\o_H(h,j_Hh) \neq 0$. This establishes the reality condition
since $ \gh = \{ A\in \gh^{\bC}| [A,j_E] = 0\}$.
\qed

In fact any tensor $S\in S^4E$ satisfying the conditions of the above lemma
can be used to define a \hkss as the following theorem shows. We 
can identify $S^4E$ with the space $\bC [E]^{(4)}$ of homogeneous 
quartic polynomials on $E \cong E^*$. 

\bt \label{t1} 
Let $S\in S^4E$, $E = \bC^{2n}$, 
be a quartic polynomial invariant under all endomorphisms
$S_{e,e'}\in S^2E = {\rm sp}(E)$ and satisfying the reality condition   

\be [S_{j_Ee,e'}-S_{e,j_Ee'},j_E] = 0\, .\la{realEqu} \ee 
Then it defines a hyper-K\"ahler
symmetric space, which is 
associated with the following complex symmetric decomposition 
\be \gg^{\bC} = \gh^{\bC} + H\ot E \, ,\quad \gh^{\bC} =  {\rm span}
\{ S_{e,e'}| e,e'\in E\} \subset {\rm sp}(E)\, .
\label{cxsymdecompEq} \ee 
The bracket $\wedge^2(H\ot E) \rightarrow \gh^{\bC}$ is given by 
(\ref{bracketEq}). The real symmetric decomposition is defined as
$\r$-real form $\gg = \gh + \gm$ of (\ref{cxsymdecompEq}), where 
\[ \gk = \gh = \{ A\in \gh^{\bC}| [A,j_E] =0\} 
= {\rm span}\{ S_{j_Ee,e'} - S_{e,j_Ee'}| e,e'\in E\}\, ,\quad 
\gm = (H\ot E)^\r.
\]  
The \hkss $M$ associated to this symmetric decomposition is the quotient
$M = M_S= G/K$, where $G$ is the simply connected Lie group with Lie algebra
$\gg$ and $K \subset G$ is the connected (and closed) subgroup with Lie 
algebra $\gk = \gh$. 

Moreover any simply connected 
\hkss can be obtained by this construction. Two hyper-K\"ahler 
symmetric spaces $M_S$ and $M_{S'}$ defined by quartics 
$S$ and $S'$ are isomorphic if and only if $S$ and $S'$ are 
in the same orbit  of the group ${\rm Aut}(E, \o_E, j_E) 
= \{ A \in {\rm Sp}(E)| [A, j_E] = 0\} \cong {\rm Sp}(k,l)$. 
\et

\pf 
First of all we note that $\gh^{\bC} =  S_{E,E} := {\rm span}
\{ S_{e,e'}| e,e'\in E\}$ is a subalgebra of ${\rm sp}(E)$ because 
\[ [S_{e,e'},S_{f,f'}] = (S_{e,e'}\cdot S)_{f,f'} - S_{S_{e,e'}f,f'}
- S_{f,S_{e,e'}f'} = - S_{S_{e,e'}f,f'}
- S_{f,S_{e,e'}f'} \in \gh^{\bC}\, .\]
Since $S$ is  $\gh^{\bC}$-invariant and completely symmetric we can check,  
 as in Lemma \ref{lemma1},  that the Jacobi identity is satisfied and 
that (\ref{cxsymdecompEq}) defines a complex symmetric decomposition.
We prove that $\gh := {\rm span}\{ S_{j_Ee,e'} - S_{e,j_Ee'}| e,e'\in E\}
\subset \{ A\in \gh^{\bC}| [A,j_E] =0\}$ defines a real form of $\gh^{\bC}$. 
Indeed for $e,e'\in E$ we have 
\[ S_{e,e'} = \frac{1}{2}(S_{e,e'}+S_{j_Ee,j_Ee'})
-\frac{\sqrt{-1}}{2}(\sqrt{-1}S_{e,e'}-\sqrt{-1}S_{j_Ee,j_Ee'}) \]
\[ = \frac{1}{2}(S_{j_Ee'',e'}-S_{e'',j_Ee'})
-\frac{\sqrt{-1}}{2}(S_{j_Ee'',\sqrt{-1}e'} - S_{e'',j_E\sqrt{-1}e'}) \, ,\]
where $e'' = -j_Ee$. Due to the reality condition the restriction of the
Lie bracket $[\cdot ,\cdot ] : \wedge^2\gm^{\bC} \rightarrow \gh^{\bC }$
to $\wedge^2\gm$ has values in $\gh$ and $\gg = \gh + \gm$ is a symmetric
decomposition with $[\gm ,\gm ] = \gh$. The metric $g^{\bC} = \o_H \ot \o_E$
defines a real valued scalar product $g$ of some signature $(p,q)$ on 
$\gm = (H\ot E)^\r$, which is 
invariant under the Lie algebra $\gh$. Since $[\gh ,j_E] = 0$ the holonomy
algebra $\gh \subset {\rm sp}(k,l)$, $p = 4k$, $q = 4l$. Hence this symmetric
decomposition defines a hyper-K\"ahler symmetric space. 

By Lemma \ref{lemma1} any \hkss can be obtained by this construction.
It is well known that a simply connected symmetric space $M$ of signature
$(p,q)$ is determined 
by its abstract curvature tensor $R \in S^2(\wedge^2 V)$, $V = \bR^{p,q}$,  
and two tensors $R$ and  $R'$ define isometric symmetric spaces if and 
only if they belong to the same ${\rm O}(V)$ orbit. Similarly  a simply 
connected hyper-K\"ahler symmetric space is determined up to isometry by its 
abstract curvature tensor  $R \in S^2(\wedge^2 V)$, where 
$V = \bR^{4k,4l}$ is the pseudo-Euclidean vector space 
with fixed hypercomplex structure $J_{\alpha} \in O(V)$. 
For a \hkss the complexified curvature tensor has the form 
\[ R(h\ot e , h'\ot e') = -\o_H(h,h') S_{e,e'} \, .\] 
where $S \in S^4E$ is the quartic of Lemma \ref{lemma1}. 
Two such curvature tensors define isomorphic hyper-K\"ahler symmetric spaces
if and only if they belong to the same orbit of 
${\rm Aut}(\bR^{4k,4l},J_{\alpha}) = {\rm Sp}(k,l)$. The group 
${\rm Sp}(k,l)$ acts on $V^{\bC} = H \ot E$ as 
${\rm Id}  \ot {\rm Sp}(E)^{j_E}
= {\rm Id}  \ot {\rm Aut}(E,\o_E,j_E)$. Hence 
two curvature tensors $R = -\o_H \ot S$ and $R' = - \o_H \ot S'$ are in the 
same ${\rm Sp}(k,l)$-orbit if and only if $S$ and $S'$ are in the same 
${\rm Sp}(k,l)$-orbit on $S^4E$. 
\qed

\section{Complex hyper-K\"ahler symmetric spaces} 
\subsection{Complex hyper-K\"ahler manifolds}
A {\bf complex Riemannian manifold} is a complex manifold $M$ equipped with
a complex metric $g$, i.e.\ a holomorphic section  $g \in \Gamma (S^2T^*M)$
which defines a  nondegenerate complex quadratic form. As in the real case
any such manifold has a unique holomorphic torsionfree and metric connection 
(Levi-Civit\`a connection). A 
{\bf complex hyper-K\"ahler manifold} is a 
complex Riemannian manifold $(M^{4n},g)$ of complex dimension $4n$ 
together with a compatible hypercomplex structure, i.e.\  three 
$g$-orthogonal parallel complex linear endomorphisms $(J_1,J_2,J_3 = J_1J_2)$ 
with $J_{\alpha}^2 = -1$.  This means that the holonomy group 
$\Hol \subset {\rm Sp}(n,\bC ) = Z_{O(4n,\bC )}({\rm Sp(1,\bC )})$.  
The linear group 
${\rm Sp}(n,\bC )$ is diagonally embedded into 
${\rm Sp}(n,\bC ) \times {\rm Sp}(n,\bC ) \subset {\rm GL}(4n, \bC)$.  
Two complex \hk manifolds $(M,g,J_{\alpha})$ 
($\alpha  = 1,2,3$)
and $(M',g',J_{\alpha}')$ are called {\bf isomorphic} if there exists a
holomorphic isometry $\varphi :M \rightarrow M'$ such that 
$\varphi^*J_{\alpha}' = J_{\alpha}$ and $\varphi^*g'= g$.

We will show that the complex \hk structure can be described as a half-flat
Grassmann structure of a certain type. A {\bf Grassmann structure} 
on a complex Riemannian manifold $(M,g)$ is a decomposition of the 
(holomorphic) tangent bundle
$TM \cong H \ot E$ into the tensor product of two holomorphic vector bundles
$H$ and $E$ of rank $2m$ and $2n$ with holomorphic nondegenerate 2-forms 
$\o_H$ and $\o_E$ such that
$g = \o_H \ot \o_E$. The Grassmann structure will be called {\bf parallel}
if the Levi-Civit\`a connection $\n = \n^{TM}$ can be decomposed as:
\[ \n = \n^H \ot {\rm Id} + {\rm Id}\ot \n^E \, ,\]
where $\n^H$ and $\n^E$ are (uniquely defined) 
symplectic connections in the bundles $H$ and $E$. 
A parallel Grassmann structure will be called {\bf half-flat} if $\n^H$ is
flat. Note that a parallel Grassmann structure on a simply connected manifold
is half-flat if and only if
the holonomy group of the Levi-Civit\`a connection is contained in 
${\rm Id}\otimes {\rm Sp}(n,\bC ) \subset {\rm Sp}(m,\bC) \ot {\rm Sp}(n,\bC)
\subset {\rm O}(\bC^{2m} \ot \bC^{2n})$.

\bp A complex hyper-K\"ahler structure $(g,J_{\alpha})$ on a simply 
connected complex manifold
$M$ is equivalent to the following geometric data:
\begin{enumerate}
\item[(i)] a half-flat Grassmann structure $(TM,g,\n )  \cong (H,\o_H,\n^H)
 \ot (E,\o_E,\n^E)$ and  
\item[(ii)] an isomorphism of flat symplectic vector bundles $H \cong M \times 
\bC^2$. Under this isomorphism $\o_H = h_1^*\wedge h_2^*$, where 
$(h_1,h_2)$ is the standard basis of $\bC^2$ 
considered as parallel frame of the trivial bundle $H = M \times 
\bC^2$. 
\end{enumerate} 
More precisely, 
\be J_1 = R_i \otimes {\rm Id}\, ,\quad J_2 = R_j\otimes {\rm Id}\, ,\quad
\mbox{and}\quad J_3 = R_k\otimes {\rm Id}\, ,\label{JalphaEq}
\ee 
where we have identified $\bC^2 = 
\bC h_1 \oplus \bC h_2$ with $\bH = {\rm span}_{\bR}\{ 1,i,j,k\} =
{\rm span}_{\bC} \{ 1,j\} = \bC 1 \oplus \bC j$ with the complex structure 
defined by  left-multiplication
by $i$ and $R_x$ denotes the right-multiplication by the quaternion $x\in \bH$.
\ep

\pf 
It is easy to check that the geometric data (i) and (ii) define a
complex hyper-K\"ahler structure on $M$. Conversely let $(g,J_{\alpha})$
be a complex hyper-K\"ahler structure on $M$. The endomorphism $J_1$ has
eigenvalues $\pm i$ and the tangent space can be decomposed into
a sum of eigenspaces 
\[ TM = E_+\oplus E_- \, .\] 
{}From the $J_1$-invariance of the metric $g$ it follows that 
$g(E_\pm ,E_\pm ) = 0$ and we can identify $E_- = E^*$ with the dual
space of $E = E_+$. Since $J_2$ anticommutes with $J_1$ it interchanges
$E$ and $E^*$ and hence defines an isomorphism 
$E \stackrel{\sim}{\rightarrow} E^*$.  Now 
$g(\cdot ,J_2 \cdot )$ defines a symplectic form $\o_E$ on $E$. 
Let $H = M \times \bC^2 = M\times (\bC h_1 \oplus \bC h_2)$ be the 
trivial bundle with 2-form $\o_H = h_1^* \wedge h_2^*$. Then we can 
identify  
\[ TM = E \oplus E^* = E \oplus E = h_1 \ot E \oplus h_2\ot  E = H \ot E
\, .\] 
We check that under this identification we have $g = \o_H \ot \o_E$. Note that
both sides vanish on $h_1 \ot E$ and $h_2\ot E$ and $\o_H(h_1,h_2) =1$. 
We calculate for $e,e' \in E = E_+ = h_1\ot E$:  
\[ g(e , J_2e' ) = \o_E(e,e') = \o_H(h_1,h_2)\o_E(e,e') = 
(\o_H \ot \o_E) (h_1\ot e, h_2\ot e')\, .\]   
Hence we have a Grassmann structure. The eigenspaces $E_\pm$ of the 
parallel endomorphism $J_1$ are invariant under parallel transport.
Therefore the Levi-Civit\`a connection $\n$ induces a connection 
$\n^E$ in the bundle
$E$. Since $\n g = 0$ and $\n J_2 = 0$ we have $\n^E \o_E = 0$. We define 
a flat connection $\n^H$ on the trivial bundle 
$H = M \times \bC^2$ by the condition $\n^H h_1 = \n^H h_2 = 0$. 
Then $\n = \n^H \ot {\rm Id} + {\rm Id} \ot \n^E$. So the Grassmann
structure is half-flat.  
 
Finally, using the standard identification $\bC^2 = \bH$,
one can easily check that the 
$J_{\alpha}$ are given by (\ref{JalphaEq}). 
\qed

\subsection{Complexification of real \hk manifolds}  
Let $(M,g,J_{\alpha})$ be a (real) hyper-K\"ahler manifold. We will 
assume that it is real analytic. This is automatically true if the metric 
$g$ is positive definite
since it is Ricci-flat and a fortiori Einstein. Using analytic continuation 
we can extend $(M,g,J_{\alpha})$ to a complex hyper-K\"ahler manifold
$(M^\bC,g^\bC,J_{\alpha}^\bC)$ equipped with an antiholomorphic involution
$T$. In complex local coordinates $z^j = x^j + iy^j$ which are extension
of real analytic coordinates $x^j, y^j$ the involution is given by the 
complex conjugation $z^j \rightarrow \bar{z}^j = x^j - iy^j$. 
We can reconstruct the (real) hyper-K\"ahler manifold as the fixed point set
of $T$. We will call $(M,g,J_{\alpha})$ a {\bf real form} of 
$(M^\bC,g^\bC,J_{\alpha}^\bC)$ and $(M^\bC,g^\bC,J_{\alpha}^\bC)$ the 
complexification of $(M,g,J_{\alpha})$. 

In general a complex hyper-K\"ahler manifold has no real form. A necessary 
condition is that the holonomy group of $\n^E$ is contained in 
${\rm Sp}(k,l)$, $n=k+l$, and hence preserves a quaternionic structure. 
Then we can define a parallel antilinear endomorphism field 
$j_E:E \rightarrow E$ such that $j_E^2 = -1$ and $\o_E (j_Ex ,j_Ey )
= \overline{\o_E (x,y)}$ for all $x,y\in E$, where the bar denotes
complex conjugation. 
We define a parallel antilinear endomorphism field $j_H: H 
\rightarrow H$  as the left-multiplication by the quaternion $j$ on 
$H = M \times \bH$. Then $\r = j_H \ot j_E$ defines a field of real 
structures in $TM = H\ot E$. We denote by ${\cal D}
\subset TM$ the real eigenspace distribution of $\r$ with eigenvalue $1$.
Here $TM$ is considered as real tangent bundle of the real manifold $M$.  
If $M^\r \subset M$ is a leaf of $\cal D$ of real dimension $4n$ then
the data $(g,J_{\alpha})$ induce on $M^\r$ a (real) \hk structure.

\subsection{Complex hyper-K\"ahler symmetric spaces} 
A complex Riemannian symmetric space is a complex Riemannian manifold 
$(M,g)$ such that any point is an isolated fixed point of an isometric 
holomorphic involution. Like in the real case one can prove that it admits
a transitive complex Lie group of holomorphic isometries and that
any simply connected complex Riemannian symmetric $M$ 
is associated to a  complex 
symmetric decomposition   
\be \gg = \gk + \gm \, ,\quad [\gk ,\gk ] \subset \gk\, ,\quad 
[\gk ,\gm ] \subset \gm \, , \quad 
[\gm ,\gm ] = \gk \label{cxsymdecomp1Eq}\ee
of a complex Lie algebra $\gg$ together with an ${\rm ad}_{\gk}$-invariant 
complex scalar product on $\gm$. More precisely  
$M = G/K$, where $G$ is the simply connected complex Lie group 
with the Lie algebra $\gg$ and $K$ is the (closed) connected subgroup 
associated with $\gk$. The holonomy group of such manifold 
is $H = {\rm Ad}_K|\gm$. Any pseudo-Riemannian 
symmetric space $M = G/K$ associated with a symmetric decomposition 
$\gg = \gk + \gm$ has a canonical complexification $M^{\bC} = G^{\bC}/K^{\bC}$
defined by the complexification $\gg^{\bC} = \gk^{\bC} + \gm^{\bC}$ of the
symmetric decomposition. Proposition \ref{fullisotropyProp} remains true for 
complex  Riemannian symmetric spaces. Ignoring the reality condition we  
obtain the following complex version of Theorem \ref{t1}. 

\bt  \la{cxt1} 
Let $S\in S^4E$, $E = \bC^{2n}$, be a quartic polynomial invariant under 
all endomorphisms
$S_{e,e'}\in S^2E = {\rm sp}(E)$. Then it defines a complex hyper-K\"ahler
symmetric space, which is 
associated with the following complex symmetric decomposition 
\be \gg = \gh + H\ot E \, ,\quad \gh =  S_{E,E} = {\rm span}
\{ S_{e,e'}| e,e'\in E\} \subset {\rm sp}(E)\, .
\label{cxsymdecomp2Eq} \ee 
The bracket $\wedge^2(H\ot E) \rightarrow \gh$ is given by 
(\ref{bracketEq}). 
The complex \hkss $M$ associated to this symmetric decomposition is 
the quotient
$M = M_S= G/K$, where $G$ is the (complex) simply connected Lie group with 
Lie algebra
$\gg$ and $K \subset G$ is the connected (and closed) subgroup with Lie 
algebra $\gk = \gh$. 

Moreover any simply connected complex \hkss can be obtained by this 
construction. 
Two complex hyper-K\"ahler 
symmetric spaces $M_S$ and $M_{S'}$ defined by quartics 
$S$ and $S'$ are isomorphic if and only if $S$ and $S'$ are in the same orbit 
of ${\rm Aut}(E, \o_E) 
= {\rm Sp}(E) \cong {\rm Sp}(n,\bC )$. 
\et

\bc \label{c1} 
There is a natural bijection between simply connected complex 
hyper-K\"ahler symmetric spaces of dimension $4n$ 
up to isomorphism and ${\rm Sp}(n,\bC )$-orbits on the space of 
quartic polynomials $S \in S^4E$ in the symplectic vector space
$E = \bC^{2n}$ such that 
\be S_{e,e'}\cdot S = 0 \quad \mbox{for all}\quad  e,e' \in E\, . 
\label{sEq} \ee
\ec 

\subsection{Classification of complex hyper-K\"ahler symmetric spaces} 
The following complex version  
of Proposition \ref{solvableProp} (with similar proof) will be  
a crucial step in the classification of complex hyper-K\"ahler symmetric 
spaces. 
\bp \label{cxsolvableProp} 
Let $(M = G/K,g,J_\a)$ be a simply connected complex 
hyper-K\"ahler symmetric space. Then the holonomy group of $M$ is solvable
and $M$ admits a transitive solvable Lie group of automorphisms. 
\ep

Due to Corollary \ref{c1} the classification of simply connected complex 
hyper-K\"ahler symmetric spaces reduces to the determination of 
quartic polynomials $S$ satisfying (\ref{sEq}). 
Below we will determine all such polynomials.  
We will prove that the following example gives all such polynomials. 

\noindent 
{\bf Example 1:} Let $E = E_+ \oplus E_-$ be a Lagrangian decomposition, 
i.e.\ $\o (E_{\pm} , E_{\pm}) = 0$, of the symplectic vector space
$E = \bC^{2n}$. Then any polynomial $S \in S^4E_+ \subset S^4E$
satisfies the condition (\ref{sEq}) and defines a simply connected complex 
hyper-K\"ahler symmetric space $M_S$ with Abelian holonomy algebra $\gh
= S_{E_+,E_+} \subset S^2E_+ \subset S^2E = {\rm sp}(E)$. 

In fact, since $E_+$ is Lagrangian the endomorphisms from $S^2E_+$ 
form an Abelian subalgebra of ${\rm sp}(E)$, which acts 
trivially on $E_+$ and hence on $S^4E_+$.

\bt \la{mainThm} 
Let $S\in S^4E$ be a quartic polynomial satisfying (\ref{sEq}). Then there 
exists a Lagrangian decomposition $E = E_+ \oplus E_-$ such that 
$S \in S^4E_+$. 
\et 
 
\pf According to Theorem \ref{cxt1} the quartic $S$ defines a \hkss with 
holonomy Lie algebra $\gh = S_{E,E}$. 
Since, by Proposition \ref{cxsolvableProp},  $\gh$ is solvable, 
Lie's theorem implies the existence of a one-dimensional $\gh$-invariant 
subspace $P = \bC p \subset E$. There exists an $\o$-nondegenerate 
subspace $W \subset E$ such that the $\o$-orthogonal complement
of $P$ is $P^{\perp} = P \oplus W$. We choose a vector $q\in E$ such that
$\o (p,q) = 1$ and $\o (W,q) = 0$ and put $Q := \bC q$. Then we have
\[ E = P \oplus  W \oplus Q \, .\] 
Since $\gh$ preserves $P$ we have the following inclusion
\[ \gh \subset PE + W^2 = P^2 + PW + PQ + W^2\, ,\] 
where we use the notation $XY = X \vee Y$ for the symmetric product of
subspaces $X,Y \subset E$. Then the second prolongation $\gh^{(2)} = 
\{ T \in S^4E| T_{e,e'} \in \gh$ for all $e,e'\in E\}$ has the following
inclusion
\be \gh^{(2)} \subset P^3E + P^2W^2 + PW^3 + W^4 = 
P^4 + P^3Q + P^3W + P^2W^2 + PW^3 + W^4  \, .\label{h2Eq} 
\ee  
Indeed $\gh^{(2)}   \subset \gh^2 = P^4 + P^3Q + P^3W + P^2Q^2 + P^2WQ +
P^2W^2 + PQW^2 + PW^3 + W^4$. The projection $\gh^{(2)} \rightarrow
P^2Q^2 + P^2WQ + PQW^2$ is zero because otherwise $S_{q,q}\in \gh 
\subset PE + W^2$ would have a nonzero projection to 
$Q^2 + WQ$ or $S_{w,q} \in \gh$ would have a nonzero projection to 
$QW$ for appropriate choice of $w\in W$. 
By (\ref{h2Eq}) we can write the quartic $S$ as
\[ S = p^3(\lambda p + \mu q + w_0) + p^2B + pC + D\, ,\]
where $\l , \mu \in \bC$, $w_0 \in W$, $B \in S^2W$, $C\in S^3W$
and  $D\in S^4W$. From now on we will identify $S^dE$ with the 
space $\bC [E^*]^{(d)}$ of homogeneous polynomials on $E^*$ of degree $d$.
Then the $\o$-contraction $T_x = \i_{\o x}T = T(\o x,\ldots )$ of a tensor 
$T \in S^dE$ with a vector
$x \in E$ is identified with the following homogeneous polynomial
of degree $d-1$:
\[ T_x = \frac{1}{d}\partial_{\omega x}T \, ,\]
where $\partial_{\omega x} T$ is the derivative of the polynomial
$T \in \bC [E^*]^{(d)}$ in the direction of $\o x = \o (x, \cdot ) \in E^*$.
For example $p_q = \langle p, \o q\rangle = \omega (q,p) = 
\partial_{\o q} p = - \partial_{p^*} p = -1 = - q_p$.

{}From $S_{p,q} = - \frac{1}{4}\mu p^2$ and 
the condition $S_{p,q} \cdot S = 0$ we obtain $\mu = 0$, since 
$p^2 \cdot S = \mu p^4$. This implies $S_{p,\cdot } = 0$.
Next we compute: 
\begin{eqnarray*}  S_{q,q} &=& \frac{1}{6}(6\lambda p^2 + 3 pw_0 + B)\\
S_{q,w} &=& -\frac{1}{12}(-3p^2\o (w_0,w) + 2p\partial_{\o w}B + 
\partial_{\o w}C)\\
& = & -\frac{1}{12}(-3p^2\o (w_0,w) + 4pB_w + 3C_w)\\  
S_{w,w'} &=& \frac{1}{6}(p^2B_{w,w'} + 3pC_{w,w'} + 6D_{w,w'}) 
\end{eqnarray*}
for any  $w,w' \in W$. 

Now the condition (\ref{sEq}) can be written as follows: 
\begin{eqnarray*} 0 & = & 
6S_{q,q}\cdot S = (3pw_0 + B)\cdot S = (\frac{3}{2}(p\otimes w_0 +
w_0 \otimes p) + B)\cdot S\\
 & = & \frac{3}{2}(2p^3B_{w_0} +
3p^2C_{w_0} + 4pD_{w_0}) + p^3Bw_0 + p^2B\cdot B + p B \cdot C +
B \cdot D\\
& = &  -2p^3Bw_0 + \frac{9}{2}p^2C_{w_0} + p (6D_{w_0}+ B \cdot C)
+ B \cdot D\, .
\end{eqnarray*}
Note that $Bw_0 = -B_{w_0}$ and $B\cdot B = [ B, B ] = 0$. 
\begin{eqnarray*} 0 & = & -12S_{q,w}\cdot S = (4pB_w + 3C_w)\cdot S\\ 
& = & 2(p^4\o (B_w,w_0) +
2p^3B^2w -3 p^2C_{Bw} -4 pD_{Bw})\\
 & & + 3(p^3C_ww_0 + p^2C_w\cdot B
+pC_w \cdot C + C_w\cdot D)\\
& = &  2p^4\o (B_w,w_0) + p^3(4B^2w + 3C_ww_0) + p^2(-6C_{Bw} + 3C_w\cdot B)\\ 
& & + p(-8D_{Bw} + 3C_w \cdot C) + 3C_w\cdot D\\
0 & = & 2S_{w,w'}\cdot S =  \frac{1}{2}(p\otimes C_{w,w'} + C_{w,w'}
\otimes p) \cdot S + 2D_{w,w'}\cdot S\\ 
& = & \frac{1}{2}(p^4\o (C_{w,w'},w_0) -2 p^3 BC_{w,w'} + 3p^2C_{C_{w,w'}} +
4pD_{C_{w,w'}}) + \\
& & 2(p^3D_{w,w'}w_0 + p^2D_{w,w'}\cdot B + pD_{w,w'}\cdot C + D_{w,w'}\cdot D)\\
& = & \frac{1}{2}p^4\o (C_{w,w'},w_0) + p^3(-BC_{w,w'} +
2D_{w,w'}w_0) +\\
& &  p^2(\frac{3}{2}C_{C_{w,w'}} 
+ 2D_{w,w'}\cdot B) 
+ p(2D_{C_{w,w'}} + 2D_{w,w'}\cdot C) + 2D_{w,w'}\cdot D
\end{eqnarray*}
This gives the following system of equations:
\begin{eqnarray*}
&(1)&  Bw_0 = 0\\
&(2)&  C_{w_0} = 0\\
&(3)& 6D_{w_0}+ B \cdot C = 0\\
&(4)& B \cdot D = 0\\ 
&(5)& \o (B_w,w_0)  = 0\\
&(6)& 4B^2w + 3C_ww_0 = 0\\
&(7)& -2C_{Bw} + C_w\cdot B = 0\\
&(8)& -8D_{Bw} + 3C_w \cdot C = 0\\
&(9)& C_w\cdot D = 0\\
&(10)& \o (C_ww',w_0) = 0\\
&(11)& -BC_ww' +
2D_{w,w'}w_0 = 0\\
&(12)& \frac{3}{2}C_{C_ww'} 
+ 2D_{w,w'}\cdot B = 0\\
&(13)& D_{C_ww'} + D_{w,w'}\cdot C = 0\\
&(14)& D_{w,w'}\cdot D = 0
\end{eqnarray*}
Note that (5) and (10) follow from (1) and (2) and that using (2)  
equation (6) says that  the endomorphism $B$ has zero square:
\[ (6') \quad B^2 = 0\, .\] 
Eliminating $D_{w_0}$ in equations (3) and (11) we obtain:
\[ (15) \quad 0 = (B\cdot C)_ww' + 3BC_ww' = 
BC_ww' - C_{Bw}w' -C_w Bw' + 3BC_ww' = 
4BC_ww'- C_{Bw}w' -C_w Bw'\, .
\] 
We can rewrite (7) as:
\[ (7') \quad -2C_{Bw}w' + C_wBw'- BC_ww' = 0\, .\]
Eliminating $C_{Bw}w'$ in  ($7'$) and (15) we obtain:
\[ (16) \quad -3BC_ww' + C_wBw' = 0\, .\]
Since the first summand is symmetric in $w$ and $w'$ we get 
\[ (17) \quad C_wBw' = C_{w'}Bw = C_{Bw}w'\, .\]
Now using  (17) we can rewrite (15) as: 
\[ (15') \quad 2BC_ww' - C_w Bw' = 0\, .\]
The equations ($15'$) and (16) show that $ BC_ww' = C_wBw' = C_{Bw}w'= 0$ and 
hence also $B\cdot C = 0$. This implies $D_{w_0} = 0$, by (3). 
Now we can rewrite (1-14) as:
\be Bw_0 = C_{w_0} = D_{w_0} = 0\ee 
\be B\cdot C = B\cdot D = 0\label{3.7Equ}\ee 
\be B^2 = 0\ee 
\be C_{Bw} = C_wB = BC_w = 0 \label{3.9Equ}\ee
\be -8D_{Bw} + 3C_w \cdot C = 0\label{3.10Equ}\ee
\be C_w\cdot D = 0\label{3.11Equ}\ee
\be \frac{3}{2}C_{C_ww'} 
+ 2D_{w,w'}\cdot B = 0 \label{3.12Equ}\ee
\be  D_{C_ww'} + D_{w,w'}\cdot C = 0 \label{3.13Equ}\ee 
\be  D_{w,w'}\cdot D = 0\la{3.14Equ} \ee
Now to proceed further we decompose $K := \mbox{ker}$ $B = W_0 \oplus W'$,
where $W_0 = \mbox{ker}$ $\o |K$ and $W'$ is a (nondegenerate) complement. 
Let us denote by $W_1$ a complement to $K$ in $W$ such that $\o (W',W_1) = 0$.
Then $W_0 + W_1$ is the $\o$-orthogonal complement to the $B$-invariant 
nondegenerate subspace $W'$. This shows that $BW_1 \subset (W_0 + W_1)\cap K
= W_0$. Moreover since $W_1 \cap K = 0$ the map $B: W_1 \rightarrow W_0$ is
injective and hence $\dim W_1 \le \dim W_0$. On the other hand 
$\dim W_1 \ge \dim W_0$, since $W_0$ is an isotropic subspace of the
symplectic vector space  $W_0 + W_1$. This shows that 
$B: W_1 \rightarrow W_0$ is an isomorphism. 

\bl
$C\in S^3K$ and $D\in S^4K$. 
\el 

\pf 
Since $W_0 = BW$ the equation (\ref{3.9Equ}) shows that  $C_{W_0} = 0$, 
which proves the 
first statement. From (\ref{3.10Equ}) and the identity 
\be (C_x\cdot C)_y = [C_x,C_y] - C_{C_xy}\label{ccEqu}\ee 
we obtain
\be D_{Bx,y} + D_{By,x} = \frac{3}{8}((C_x\cdot C)_y + (C_y\cdot C)_x) =
-\frac{3}{4}C_{C_xy}\, .\label{sym1Equ} \ee 
The equation $B\cdot D = 0$ (\ref{3.7Equ}) reads:
\[ 0 = (B\cdot D)_{x,y} = [B,D_{x,y}] - D_{Bx,y} -D_{x,By}\, . \]
Using this (\ref{3.12Equ}) yields:
\be D_{Bx,y} + D_{By,x} = [B,D_{x,y}] =  - D_{x,y}\cdot B 
= \frac{3}{4} C_{C_xy}\, .\label{sym2Equ}\ee
Now from (\ref{sym1Equ}) and (\ref{sym2Equ}) we obtain that
\[ 0 = C_{C_xy}z = C_zC_xy \] 
for all $x,y,z \in W$. This implies $[C_x,C_y]  = 0$ for all $x,y\in W$ 
and hence
\be C_x \cdot C = 0\ee
for all $x\in W$, by (\ref{ccEqu}). Finally this shows that 
$D_{W_0} = 0$ by (\ref{3.10Equ}).  This proves the second statement.
\qed 

\bl $D_{x,y}C_z = C_zD_{x,y} = 0$ for all $x,y,z\in W$. 
\el 

\pf Using (\ref{3.13Equ}) we compute:
\be D_{x,y}C_z w = D_{C_zw,x}y = -(D_{z,w}\cdot C)_xy = 
-([D_{z,w},C_x]y - C_{D_{z,w}x}y) \, .\label{dcEqu}\ee 
{}From (\ref{3.11Equ}) we get:
\[ 0 = (C_x\cdot D)_{z,w} y = [C_x,D_{z,w}]y - D_{C_xz,w}y - D_{z,C_xw}y =
C_xD_{z,w}y - D_{z,w}C_xy -  D_{y,w}C_xz - D_{z,y}C_xw \, ,\]
and hence:
\[ [D_{z,w},C_x]y = -  D_{y,w}C_xz - D_{z,y}C_xw\, ,\]
and 
\[ C_{D_{z,w}x}y = C_yD_{z,w}x =  D_{z,w}C_xy +  
D_{x,w}C_yz + D_{z,x}C_yw \, .\] 
Now we eliminate the $CD$-terms from (\ref{dcEqu}) arriving at:
\be  D_{x,y}C_z w = (D_{y,w}C_xz + D_{z,y}C_xw + D_{z,w}C_xy +    
D_{x,w}C_yz + D_{z,x}C_yw) \, . \label{dxyzwEqu}\ee  
Considering all the permutations of $(x,y,z,w)$ we get 6 homogeneous  linear
equations for the 6 terms of equation  (\ref{dxyzwEqu}) with the 
matrix:
\[ \left( \begin{array}{cccccc}
-1 & 1 & 1 & 1 & 1 & 1\\
 1 & -1 & 1 & 1 & 1 & 1\\
 1 & 1 & -1 & 1 & 1 & 1\\
 1 & 1 & 1 & -1 & 1 & 1\\
 1 & 1 & 1 & 1 & -1 & 1\\
 1 & 1 & 1 & 1 & 1 & -1\\
\end{array}\right) \] 
This is the matrix of the endomorphism $-2\ {\mbox{Id}} + e\otimes e$ in 
the arithmetic space $\bR^6$, where $e = e_1 + \dots + e_6$; 
$(e_i)$ the standard basis. It has eigenvalues $(4,-2,-2,-2,-2,-2)$. 
This shows that the matrix is nondegenerate and proves the lemma.
\qed 

For a symmetric tensor $T\in S^dW$ we denote by 
\[ \Sigma_T := 
{\rm span}\{ T_{x_1,x_2,\ldots , x_{d-2}}x_{d-1}| x_1,x_2, \ldots 
x_{d-1} \in W\} \subset W\}\]
the {\bf support} of $T$. 

\bl \la{supportLemma} 
The supports of the tensors $B\in S^2W$, $C\in S^3W$ and $D\in S^4W$ 
admit the following inclusions 
\[ \Sigma_B + \Sigma_C \subset 
{\rm ker} B\cap {\rm ker} C \cap {\rm ker} D\, ,
\quad \Sigma_D  \subset {\rm ker} B\cap {\rm ker} C\, .\]
Moreover $\Sigma_B + \Sigma_C$ is isotropic and $\o (\Sigma_D,
\Sigma_B + \Sigma_C) = 0$. 
\el

\pf 
The first statement follows from $B^2 = BC_x = BD_{x,y}= C_xB = C_xC_y = 
C_xD_{y,z} = D_{x,y}B = D_{x,y}C_z = 0$ for all $x,y,z\in W$. The second 
statement follows from the first and the definition of support, e.g.\
if $z = C_xy\in \Sigma_C$ and $w\in \Sigma_B + \Sigma_C + \Sigma_D 
\subset {\rm ker} C$ we compute:
\[ \o (z,w) = \o (C_xy,w) = -\o (y,C_xw) = 0\, .\]
\qed

\bl \la{solvableLemma} 
The Lie algebra $D_{W,W} \subset S^2W \cong {\rm sp}(W)$ is solvable.
\el

\pf This follows from Proposition \ref{cxsolvableProp}, since $D \in S^4W$
satisfies (\ref{sEq}) and hence defines a complex hyper-K\"ahler symmetric
space with holonomy Lie algebra $D_{W,W}$. It also follows from the 
solvability of $S_{E,E}$ as we show now.  
In terms of the decomposition $E = P + W +Q$ an endomorphism 
\[ S_{x,y} = (\lambda p^4 + p^3w_0 + p^2B + pC + D)_{x,y} = 
(p^2B + pC + D)_{x,y} = B(x,y)p^2 - pC_xy + D_{x,y}\, ,\] 
where $x, y\in W$,  is represented by 
\[ \left( \begin{array}{ccc} 0& -(C_xy)^t & B(x,y)\\
                             0& D_{x,y} & -C_xy\\
                             0& 0       & 0
\end{array}
\right)\, .\] 
Since the Lie algebra $S_{E,E}$ is solvable this implies 
that the Lie algebra  $D_{W,W}$, which corresponds to the induced
representation of $S_{E,E}$ on $P^{\perp}/P \cong W$, is also solvable. 
\qed 

\bl \la{w0Lemma}
\[ \o (w_0, \Sigma_B + \Sigma_C + \Sigma_D ) = 0\, .\]
\el

\pf 
Note that $w_0 \in {\rm ker} B \cap {\rm ker} C \cap {\rm ker} D$, due to
equations (1-3) and (\ref{3.7Equ}). This implies the lemma. In fact, if 
e.g.\ $y = Bx \in \Sigma_B$ then
\[ \o (w_0, y) = \o (w_0,Bx) = -\o (Bw_0,x) = 0\, ,\]
which shows that $\o (w_0,\Sigma_B) = 0$ 
\qed

Now to finish the proof of Theorem \ref{mainThm} we will use induction on the
dimension $\dim E = 2n$. If $n=1$ the (solvable) holonomy algebra $\gh$ 
is a proper
subalgebra of $S^2E \cong {\rm sl}(2,\bC)$. Without loss of generality we 
may assume that either\\
a) $\gh = \bC p^2$ or\\
b) $\gh = \bC pq$ or\\ 
c) $\gh = \bC p^2 +  \bC pq$,\\
where $(p,q)$ is a symplectic basis of $E$. 
In the all three cases the Lie algebra $S_{E,E} = \gh \subset 
\bC p^2 +  \bC pq$ and
hence
\[ S = \l p^4 + \mu p^3q\, .\]
In the cases b) and c) we have that $pq \in \gh$ and since  
\[ pq \cdot S = \frac{1}{2}(p \ot q + q\ot p) \cdot S 
= -2\l p^4 -\mu p^3q = 0\]
it follows that $S=0$. In the case a) from $S_{E,E}  = \gh = \bC p^2$ we
have that $S = \l p^4$. This tensor is invariant under  $\gh = \bC p^2$
and belongs to the fourth symmetric power of the Lagrangian
subspace $\bC p \subset E$. This establishes the first step of
the induction. Now by induction using 
equation (\ref{3.14Equ}) and Lemma \ref{solvableLemma} 
we may assume that $\Sigma_D$ is isotropic. Now Lemma \ref{supportLemma}
and Lemma \ref{w0Lemma} show that $\bC w_0 + \Sigma_B + \Sigma_C + \Sigma_D$
is isotropic and hence is contained in some Lagrangian subspace $E_+ \subset 
E$. This implies that $S\in S^4E_+$. 
\qed

Now we give a necessary and sufficient condition for a symmetric
manifold $M = M_S$, $S \in S^4E_+$, to have no flat de Rham factor.

\bp  \la{flatProp} 
The complex \hkss $M_S$, $S \in S^4E_+$, has no flat de Rham factor if 
and only if the support $\Sigma_S = E_+$. 
\ep 

\pf 
If $M = M_S =G/K$ has a flat factor $M_0$, such that $M = M_1 \times M_0$, 
then this induces a decomposition $E = E^1 \oplus E^0$ and 
$S \in S^4E_1$; hence $\Sigma_S \subset E_1\cap E_+ \neq E_+$. Conversely let
$S\in S^4 E_+$, assume that $E^1_+ = \Sigma_S \subset E_+$ is a proper 
subspace and choose a complementary subspace $E^0_+$. We denote by
$E_-^1$ and $E_-^0$ the annihilator of $\o E_+^0$ and $\o E_+^1$
respectively. Let us denote $E^1 = E^1_+ \oplus E^1_-$, 
$E^0 =  E^0_+ \oplus E^0_-$, $\gm^1 = H \otimes E^1$ and 
$\gm^0 = H \otimes E^0$. Then $E^0, E^1 \subset E$ are $\o$-nondegenerate
complementary 
subspaces and $\gm^0,\gm^1 \subset \gm = T_oM$ are $g$-nondegenerate 
complementary subspaces. Since $S \in S^4E_+^1$ the Lie algebra 
$\gg = \gh + \gm = (\gh + \gm^1) \oplus \gm^0$ has 
the Abelian direct 
summand $\gm^0$, see (\ref{bracketEq}), which gives rise to a flat factor $M^0 \subset M_S
= M^1 \times M^0$. 

\bt \la{cxmainThm} 
Any simply connected complex hyper-K\"ahler symmetric space 
without flat de Rham factor is isomorphic to a complex hyper-K\"ahler 
symmetric 
space of the form $M_S$, where 
$S \in S^4E_+$ and $E_+ \subset E$ is a 
Lagrangian subspace of the complex symplectic vector space $E = \bC^{2n}$. 
Moreover there is a natural 
1-1 correspondence between simply connected complex hyper-K\"ahler symmetric 
spaces  
without flat factor up to isomorphism 
and orbits $\cal O$ 
of the group ${\rm Aut}(E,\o , E_+) = \{ A\in Sp(E)|AE_+ = E_+\} \cong 
{\rm GL}(E_+) \cong {\rm GL}(n,\bC )$ on the space $S^4E_+$ 
such that $\Sigma_S = E_+$ for all $S\in {\cal O}$. 
\et 

\pf This is a corollary of Theorem \ref{cxt1}, Theorem \ref{mainThm} and
Proposition  \ref{flatProp}. \qed

Let $M = G/K$ be a simply connected complex hyper-K\"ahler symmetric 
space without flat factor. By Theorem \ref{mainThm}
and Proposition \ref{flatProp} it is associated to quartic polynomial 
$S\in S^4E_+$ with support $\Sigma_S = E_+$. Now we describe the Lie algebra 
${\rm aut} (M_S)$ of the full  group of automorphisms, i.e.\ 
isometries which preserve the hypercomplex structure, of $M_S$.

\bt \la{fullautThm} Let $M_S = G/K$ be as above. Then the full automorphism 
algebra is given by
\[ {\rm aut}(M_S) = {\rm aut}(S)  + \gg \, ,\]
where $A \in {\rm aut}(S) = \{ B\in {\rm gl}(E_+)| B \cdot S = 0\}$ acts on 
$\gg = \gh + \gm$ as follows. It preserves the decomposition and acts on 
$\gh = S_{E,E}$ by 
\[ [A,S_{x,y}] = S_{Ax,y} + S_{x,Ay} \]
for all $x,y\in E$ and on $\gm = H\ot E$ by 
\[ [A, h\ot e] = h \ot Ae \, ,\]
where ${\rm gl}(E_+)$ is canonically embedded into ${\rm sp}(E)$.  
\et 

\pf 
By (the complex version of) 
Proposition \ref{fullisotropyProp} it is sufficient to determine the 
centralizer $\gc$ of ${\rm sp}(1,\bC)$ in the full isotropy algebra 
$\tilde{\gh} 
= {\rm aut}(R) \supset {\rm sp}(1, \bC) \oplus \gh$. Equation
(\ref{fullisotropyEq}) shows that 
\[ \gc =  \{ {\rm Id} \ot A|\,   A\in {\rm sp}(E)\, ,\;  
A \cdot S = [A, S(\cdot ,\cdot )] 
- S(A \cdot ,\cdot ) 
 - S(\cdot ,A\cdot ) = 0\}  \,. \]
{}From $A \cdot S = 0$ we obtain that the commutator   
$[A,S_{x,y}] = S_{Ax,y} + S_{x,Ay}$ for all $x,y\in E$ and 
$A \Sigma_S = AE_+ \subset E_+$. This implies  $\gc = {\rm aut}(S)$. 
\qed

\section{Classification of hyper-K\"ahler symmetric spaces}
Using the description of complex hyper-K\"ahler symmetric spaces
given in Theorem \ref{cxmainThm} we will now
classify (real) hyper-K\"ahler symmetric spaces. Recall that a simply connected pseudo-Riemannian 
manifold is called {\bf indecomposable} if it is not a Riemannian product
of two pseudo-Riemannian manifolds. Any simply connected pseudo-Riemannian 
manifold can be decomposed into the Riemannian product of indecomposable
pseudo-Riemannian manifolds. By Wu's theorem \cite{W} a 
simply connected pseudo-Riemannian manifold is indecomposable if and only 
if its holonomy group is {\bf weakly irreducible}, i.e.\ has no 
invariant proper nondegenerate subspaces.  
Therefore it is sufficient to classify 
(real) hyper-K\"ahler symmetric spaces
with {\em indecomposable} holonomy.

Let $(M = G/K,g,J_{\alpha})$ be a hyper-K\"ahler symmetric space   
associated to a symmetric decomposition
(\ref{symdecomp}). The complexified tangent space of $M$ is identified
with $\gm^{\bC} = H \ot E$, the tensor product of to complex symplectic
vector spaces with quaternionic structure $j_H$ and $j_E$ such that
$\r = j_H\ot j_E$ is the complex conjugation of $\gm^{\bC}$ with
respect to $\gm$. By Theorem \ref{t1} it is defined by a quartic 
polynomial $S \in S^4E$ satisfying the conditions of the
theorem. Moreover the 
holonomy algebra $\gh$ acts trivially on $H$ and is identified with 
the real form of the complex Lie algebra
$S_{E,E}\subset {\rm sp}(E)$ given by 
$\gh = {\rm span}\{ S_{je,e'} -S_{e,je'}|
e, e' \in E\} = \{ A \in S_{E,E}| [A,j_E] = 0\}\subset {\rm sp}(E)^{j_E}$.  

The quartic polynomial $S$ defines also a complex hyper-K\"ahler
symmetric space $M^{\bC} = G^{\bC}/K^{\bC}$, which is 
the complexification of $M = G/K$. By Theorem \ref{mainThm}, $S \in S^4L$
for some Lagrangian subspace $L\subset E$. Recall that the symplectic
form $\o = \o_E$ together with the quaternionic structure  $j = j_E$ define a 
Hermitian metric $\gamma = \gamma_E = \o_E (\cdot , j_E\cdot )$ of 
(real) signature $(4k,4l)$, $n = k+l$, which coincides with the 
signature of the pseudo-Riemannian metric $g$ (we normalize 
$\gamma_H = \o_H (\cdot , j_H\cdot )$ to be positive definite).
We may decompose $\gamma$-orthogonally $L = L^0 \oplus L^+ \oplus L^-$, 
such that
$\gamma$ vanishes on $L^0$ is positive definite on $L^+$ and 
negative definite on $L^-$. 

\bl \la{L0+-Lemma} 
\begin{enumerate}
\item[(i)] $j L^0 = L^0$ and 
\item[(ii)] $L^+ + L^- + jL^+ + jL^- \subset E$ is an $\o$-nondegenerate
and $\gamma$-nondegenerate $\gh$-invariant subspace (with trivial action
of $\gh$). 
\end{enumerate} 
\el

\pf 
We show first that $L + jL^0$ is $\o$-isotropic and hence $L + jL^0 = L$
since $L$ is Lagrangian. Indeed $L \supset L^0$ is $\o$-isotropic and
also $jL^0$ because $\o$ is $j$-invariant. So it suffices to remark 
that $\o (L,jL^0) = 0$:
\[ \o (L,jL^0) = \gamma (L,L^0) = 0\, .\]
This implies that $jL^0 \subset L$. Since 
$\gamma (L,jL^0) = -\o (L, L^0) = 0$, we conclude that  
$jL^0 \subset {\rm ker} \gamma|L  = L^0$. This proves (i).

To prove (ii) it is sufficient to check that the subspace 
$L^+ + L^- + jL^+ + jL^-\subset E$ is 
nondegenerate with respect to $\gamma$,  since it is $j$-invariant. 
First we remark that $\gamma$ is positive definite on 
$L^+$ and $jL^+$ and negative definite on $L^-$ and $jL^-$, due to the 
$j$-invariance of $\gamma$: $\g (jx,jx) = \g (x,x)$, $x\in E$. 
So to prove (ii) it is sufficient to check
that  $jL^+ \oplus jL^-$ is $\g$-orthogonal to the $\gamma$-nondegenerate
vector space $L^+ + L^-$: 
\[ \gamma (L^+ +L ^-,jL^+ + jL^- ) = \o ( L^+ + L^-,L^+ + L^-) = 0\, .\]
\qed

By Theorem \ref{t1} the quartic polynomial $S$ must satisfy the reality 
condition $[S_{je,e'}-S_{e,je'},j] = 0$. Now we describe all such polynomials.

The quaternionic structure $j$ on $E$ is {\bf compatible} with $\o$, i.e.\ 
$\o (j x ,j y ) = \overline{\o (x,y)}$ for all $x,y\in E$ and it  
induces a real structure (i.e.\ an 
antilinear involution) $\tau := j\ot j \ot \cdots \ot j$ on all even powers 
$S^{2r}E \subset E\ot E \ot \cdots \ot E$. For $S\in S^{2r}E$ and 
$x_1,\cdots ,x_{2r} \in E$ we
have
\[ (\tau S)(x_1,\cdots ,x_{2r}) = \overline{S(jx_1,\cdots ,jx_{2r})}\, .
\]
Note that the fixed point set ${\rm sp}(E)^{\tau} =
\{ A\in {\rm sp}(E)| [A,j] = 0\} \cong {\rm sp}(k,l)$. 

\bp \la{realityProp} Let $(E,\o ,j)$ be a complex symplectic vector space 
 with a quaternionic structure $j$ such that
$\o (j x ,j y ) = \overline{\o (x,y)}$ for all $x,y\in E$. Then a quartic
polynomial $S\in S^4E$ satisfies the reality condition 
$[S_{je,e'}-S_{e,je'},j] = 0$ if and only if $S \in (S^4E)^{\tau} = 
{\rm span}\{ T +\tau T| T\in S^4E\}$. 
\ep 

\pf 
The reality condition for $S\in S^4E$ can be written as
\[ [S_{jx,jy} + S_{x,y},j]z = 0 \] 
for all $x,y,z \in E$. Contracting this vector equation with $jw\in E$ 
by means of $\o$ and using the compatibility  between 
$j$ and $\o$ we obtain the equivalent condition
\begin{eqnarray} 0  &=& -\o (jw,[S_{jx,jy} + S_{x,y},j]z)\nonumber \\
&=&  S(jx,jy,jz,jw) - 
\overline{S(jx,jy,z,w)} + S(x,y,jz,jw) - \overline{S(x,y,z,w)}\, .
\la{realityEqu} 
\end{eqnarray}   
Now putting $x=y=z=w=u$ we obtain:
\[ 0 = S(ju,ju,ju,ju) - 
\overline{S(ju,ju,u,u)} + S(u,u,ju,ju) - \overline{S(u,u,u,u)}\, \]  
and putting $x= iu$ and $y=z=w=u$ we obtain:
\[ 0 = -iS(ju,ju,ju,ju) - 
i\overline{S(ju,ju,u,u)} + iS(u,u,ju,ju) + i\overline{S(u,u,u,u)}\, .\] 
Comparing these two equations we get $S(ju,ju,juj,u) = \overline{S(u,u,u,u)}$,
i.e.\ $S = \tau S$.  This shows that the reality condition implies
that $S\in (S^4E)^{\tau}$. Conversely the condition $S = \tau S$  can
be written as
\[ S(jx,jy,jz,jw) = \overline{S(x,y,z,w)} \quad \mbox{for all}
\quad  x,y,z,w \in E\, .\]
Changing $z\rightarrow jz$ and $w\rightarrow jw$ in this equation we 
obtain 
\[ S(jx,jy,z,w) = \overline{S(x,y,jz,jw)} \quad \mbox{for all}
\quad  x,y,z,w \in E\, .\]
These two equations imply (\ref{realityEqu}) and hence the reality condition.
\qed

Now we are ready to classify simply connected 
hyper-K\"ahler symmetric spaces. 
We will show that the following construction
gives all such symmetric spaces.

Let $(E,\o ,j)$ be a complex symplectic vector space 
of dimension $2n$ with a quaternionic structure $j$ such that
$\o (j x ,j y ) = \overline{\o (x,y)}$ for all $x,y\in E$ 
and $E = E_+ \oplus E_-$ a $j$-invariant 
Lagrangian decomposition. Such a decomposition exists if and only if 
the Hermitian form $\gamma = \o (\cdot , j\cdot )$  has real signature
$(4m,4m)$, where $\dim_{\bC} E = 2n = 4m$.  Then any polynomial 
$S \in (S^4E_+)^{\tau} = S^4E_+\cap (S^4E)^{\tau}$ 
satisfies the condition (\ref{sEq}) and the reality condition,  
by Proposition \ref{realityProp}. Hence by Theorem \ref{t1} it 
defines a (real) simply connected  
hyper-K\"ahler symmetric space $M_S$ with Abelian holonomy algebra $\gh
= (S_{E_+,E_+})^{\tau} = S_{E_+,E_+}\cap (S^2E)^{\tau} = 
{\rm span}\{S_{je,e'} -S_{e,je'}|
e,e'\in E\} \subset  {\rm sp}(E)^{\tau} \cong {\rm sp}(m,m)$.

\bt \la{realmainThm} 
Any simply connected hyper-K\"ahler symmetric space 
without flat de Rham factor is isomorphic to a hyper-K\"ahler symmetric 
space of the form $M_S$, where 
$S = T + \tau T$, $T \in S^4E_+$ and $E_+ \subset E$ is a $j$-invariant
Lagrangian subspace of the complex symplectic vector space $E$ with 
compatible quaternionic structure $j$. A hyper-K\"ahler symmetric 
space of the form $M_S$ has no flat factor if and only if it complexification
has no flat factor, which happens if and only if the support $\Sigma_S = E_+$.
Moreover there is a natural 
1-1 correspondence between simply connected hyper-K\"ahler symmetric spaces  
without flat factor up to isomorphism 
and orbits $\cal O$ 
of the group ${\rm Aut}(E,\o,j,E_+) = \{ A\in Sp(E)| [A,j] = 0, 
AE_+ = E_+\} \cong {\rm GL}(m,\bH )$ on the space $(S^4E_+)^{\tau}$
such that $\Sigma_S = E_+$ for all $S\in {\cal O}$. 
\et 

\pf 
Let $M$ be a simply connected hyper-K\"ahler symmetric space. 
We first assume that it is indecomposable. Then the holonomy 
algebra $\gh$ is weakly irreducible. By Theorem \ref{t1}, $M = M_S$ for some
quartic polynomial $S \in S^4E$ satisfying (\ref{sEq}) and the 
reality condition (\ref{realEqu}).  By Proposition \ref{realityProp}
the reality condition means that $S \in (S^4E)^{\tau}$. On the other
hand, by Theorem \ref{mainThm} $S \in S^4L$ for some Lagrangian
subspace $L$ of $E$. Now the weak irreducibility of $\gh$ and
Lemma~\ref{L0+-Lemma} imply that $L = L^0$ is $j$-invariant. 
This proves that $S \in (S^4E_+)^{\tau}$, where $E_+ = L = L^0$ is 
a $j$-invariant Lagrangian subspace of $E$. This shows that 
$M$ is obtained from the above construction. Any simply connected \hkss 
$M$ without flat factor is the Riemannian product of indecomposable ones, say 
$M = M_1\times M_2 \times \cdots \times M_r$, and we may assume
that $M_i = M_{S_i}$, $S_i \in S^4E_i$. Therefore $M$ is associated
to the quartic polynomial $S = S_1 \oplus S_2 \oplus \cdots \oplus S_r
\in S^4E$, $E = E_1 \oplus E_2 \oplus \cdots \oplus E_r$. 
Moreover $S$ satisfies (\ref{sEq}) and (\ref{realEqu}) if the $S_i$ 
satisfy (\ref{sEq}) and (\ref{realEqu}). This shows that any simply connected
\hkss is obtained from the above construction.

It is clear that the complexification $M_S^{\bC}$ has a flat factor
if $M_S$ has a flat factor. Conversely let us assume that $M_S^{\bC}$ has a 
flat factor, hence $\Sigma_S \subset E_+$ is a proper subspace. Since
$j\Sigma_S = jS_{E,E}E = S_{E,E}jE = S_{E,E}E = \Sigma_S$ there exists a
$j$-invariant complementary subspace $E_+'$ in $E_+$. Denote by 
$E_-'$ the annihilator of $\Sigma_S$ in $E_-$ then $E' = E_+' \oplus E_-'$
is an $\o$-nondegenerate and $j$-invariant subspace of $E$ on which the
holonomy $\gh^{\bC} = S_{E,E} \subset S^2\Sigma_S$ acts trivially. Then the 
corresponding real subspace $(H\ot E')^{\r} \subset \gm = (H\ot E)^{\r}$ 
is a $g$-nondegenerate subspace on which the holonomy $\gh$ acts trivially. 
By Wu's theorem \cite{W} it defines a flat de Rham factor. 

Now the last statement follows from the corresponding statement in 
Theorem \ref{t1}. 
\qed

\bc Any \hkss without flat factor has signature
$(4m,4m)$. In particular its dimension is divisible by $8$. 
\ec 

\bc \label{divisibleCor} Let $M = M_S$ be a complex \hkss without 
flat factor associated with a quartic $S \in S^4E_+$, where $E_+ \subset E$
is a Lagrangian subspace. It admits a real form if and only if there exists
a quaternionic structure $j$ on $E$ compatible with $\o$ preserving $E_+$
such that $\tau S = S$, where $\tau$ is the real structure on $S^4E$ induced
by $j$. In particular $\dim_{\bC} M$ has to be divisible by $8$. 
\ec

Let $M = G/K$ be a simply connected hyper-K\"ahler symmetric 
space without flat factor. By Theorem 
\ref{realmainThm} it is associated to a 
quartic polynomial $S\in (S^4E_+)^{\tau}$ with support 
$\Sigma_S = E_+$. Now we describe the Lie algebra 
${\rm aut} (M_S)$ of the full  group of automorphisms, i.e.\ 
isometries which preserve the hypercomplex structure of $M_S$.

\bt Let $M_S = G/K$ be as above. Then the full automorphism algebra is given by
\[ {\rm aut}(M_S) = {\rm aut}(S)  + \gg \, ,\]
where ${\rm aut}(S) = \{ A\in {\rm gl}(E_+)|\, [A,j] = 0\, ,
\; A \cdot S = 0\}$ acts on 
\[ \gg = \gh + \gm \, ,\quad  \gh = \{ A \in S_{E,E}| [A,j] = 0\} = 
{\rm span} \{ S_{jx,y} - S_{x,jy}| x,y\in E\}\, ,\gm = (H\ot E)^{\r}\] 
as in Theorem \ref{fullautThm}. 
\et 

\pf 
The proof is similar to that of Theorem \ref{fullautThm}.  
\qed

\section{Low dimensional hyper-K\"ahler symmetric spaces}
\subsection{Complex hyper-K\"ahler symmetric spaces of 
dimension $\le 8$} 
{\bf Dimension 4}\\ 
Assume that $M$ is a simply connected
complex hyper-K\"ahler symmetric space of dimension
4.  Applying 
Theorem \ref{cxmainThm} we conclude
that $M = M_S$ for some $S \in S^4E_+$, where $E_+ \subset E$ is a 
one-dimensional subspace $E_+ = \bC e$. This proves:

\bt There exists up to isomorphism only one non-flat simply connected 
complex hyper-K\"ahler symmetric space of dimension 4: $M = M_S$
associated with the quartic $S = e^4$. 
\et 

\noindent 
{\bf Dimension 8}\\
Any eight-dimensional simply connected 
complex hyper-K\"ahler symmetric space is associated with a 
quartic $S \in S^4E_+$, where $E_+ \subset E$
is a Lagrangian subspace of $E = \bC^4$. 
We denote by $(e,e')$ a basis of $E_+$. 

\bt Eight-dimensional simply connected 
complex hyper-K\"ahler symmetric space are in natural 
1-1 correspondence with the orbits of the group
${\rm CO}(3,\bC ) = \bC^* \cdot {\rm SO}(3,\bC )$ on the space 
$S^2_0\bC^3$ of traceless symmetric matrices. The complex hyper-K\"ahler 
symmetric space associated with a traceless symmetric matrix $A$ is the 
manifold $M_{S(A)}$, where $S(A) \in S^4\bC^2$ is the quartic polynomial
which corresponds to $A$ under the ${\rm SO}(3,\bC )$-equivariant 
isomorphism $S_0^2\bC^3 \cong$\linebreak[4] $S_0^2\wedge^2\bC^3 \cong S_0^2S^2\bC^2 = 
S^4\bC^2$.
\et 

The classification of ${\rm SO}(3,\bC )$-orbits on 
$S^2_0\bC^3$ was given by Petrov \cite{P} 
in his classification of Weyl tensors
of Lorentzian 4-manifolds. 

\pf
By Theorem \ref{cxmainThm} the classification of 
eight-dimensional simply connected complex hyper-K\"ahler symmetric spaces
reduces to the description of orbits of the group ${\rm GL}(E_+) = 
{\rm GL}(2,\bC)$
on $S^4{\bC}^2 \subset S^2S^2{\bC}^2$. Fixing a volume form $\s$ on 
${\bC}^2$ we can identify $S^2{\bC}^2$ with ${\rm sp}(1,\bC ) \cong 
{\rm so}(3,\bC )$. Then the 
Killing form $B$ is an ${\rm SL}(2,\bC)$-invariant and we have the
${\rm GL}(2,\bC)$-invariant decomposition: $S^2S^2{\bC}^2 = S^2_0S^2{\bC}^2
\oplus \bC B$. The action of ${\rm SL}(2,\bC)$ on $S^2{\bC}^2$ is 
effectively equivalent to the adjoint action of ${\rm SO}(3,\bC )$. 
The problem thus reduces essentially to the determination of the orbits of
${\rm SO}(3,\bC )$ on $S_0^2{\bC}^3$.
\qed

\subsection{Hyper-K\"ahler symmetric spaces of 
dimension $\le 8$} 
By Corollary \ref{divisibleCor} the minimal dimension of 
non-flat hyper-K\"ahler symmetric spaces is 8. 

\bt Eight-dimensional simply connected 
hyper-K\"ahler symmetric space are in natural 
1-1 correspondence with the orbits of the group
$\bR^+ \cdot {\rm SO}(3)$ on the space 
$S^2_0\bR^3$ of traceless symmetric matrices. The hyper-K\"ahler 
symmetric space associated with a traceless symmetric matrix $A$ is the 
manifold $M_{S(A)}$, where $S(A) \in (S^4\bC^2)^{\tau}$ is 
the quartic polynomial
which corresponds to $A$ under the ${\rm SO}(3)$-equivariant 
isomorphism $S_0^2\bR^3 \cong (S^4\bC^2)^{\tau}$.
\et

\end{document}